\documentclass[preprint,12pt]{elsarticle}
\usepackage{times}
\usepackage{amsfonts,amsmath,amssymb,amsthm}
\usepackage{color,soul}
\biboptions{numbers,sort&compress} 
\journal{arXiv-28-Marh-2017}

\let\ts=\textstyle
\def\fracskip{\mskip 1mu \relax}
\let\oldfrac=\frac
\def\nfrac#1#2{\oldfrac{\fracskip#1\fracskip}{\fracskip#2\fracskip}}
\def\tfrac#1#2{{\ts\nfrac{#1}{#2}}}
\let\frac=\nfrac
\def\pd#1#2{\frac{\partial#1}{\partial#2}}
\def\pdd#1#2#3{\ifx#2#3\pd{^2#1}{#2^2}\else\pd{^2#1}{#2\partial#3}\fi }
\def\dash{\unskip\nobreak\hskip.05555em---\hskip.05555em\relax}

\let\bl=\bigl \let\br=\bigr

\def\Equation#1. {\medbreak{\bfseries\itshape{Equation\kern.3333em\relax#1.}}\enspace\ignorespaces }

\newtheoremstyle{remark}{\medskipamount}{\medskipamount}
  {\small\rmfamily}{\parindent}{\footnotesize\sffamily}{.}{.5em}{}
\theoremstyle{remark}
\newtheorem{remark}{Remark}

\let\ds=\displaystyle
\let\ts=\textstyle
\let\bl=\bigl \let\br=\bigr

\begin{document}
\begin{frontmatter}
\title{New numerical methods for blow-up problems}
\author[ipm,bmstu,mephi]{Andrei D. Polyanin\corref{cor1}}
  \ead{polyanin@ipmnet.ru}
\author[us]{Inna K. Shingareva\corref{cor2}}
  \ead{inna@mat.uson.mx}
  \cortext[cor1]{Principal corresponding author}
  \cortext[cor2]{Corresponding author}
\address[ipm]{Institute for Problems in Mechanics, Russian Academy of Sciences,\\
  101 Vernadsky Avenue, bldg~1, 119526 Moscow, Russia}
\address[bmstu]{Bauman Moscow State Technical University,\\
  5 Second Baumanskaya Street, 105005 Moscow, Russia}
\address[mephi]{National Research Nuclear University MEPhI,
  31 Kashirskoe Shosse, 115409 Moscow, Russia}
\address[us]{University of Sonora, Blvd. Luis Encinas y Rosales S/N, Hermosillo C.P. 83000, Sonora, M\'exico}

\begin{abstract}
Two new methods of numerical integration of Cauchy problems for ODEs with blow-up solutions are described.
The first method is based on applying a differential transformation,
where the first derivative (given in the original equation) is chosen as a new independent variable.
The second method is based on introducing a new non-local variable that reduces ODE to
a system of coupled ODEs.
Both methods lead to problems whose solutions do not have blowing-up singular points;
therefore the standard numerical methods can be applied.
The efficiency of the proposed methods is illustrated with several test problems.
\end{abstract}

\end{frontmatter}

\section{Introduction}

We will consider Cauchy problems for ODEs, whose solutions tend to infinity at some finite value of~$x$,
say $x=x_*$. The point $x_*$ is not known in advance.
Similar solutions exist on a bounded interval and are called blow-up solutions.
This raises the important question for practice: how to determine the position of a singular point~$x_*$
and the solution in its neighborhood by numerical methods.
In general, the blow-up solutions, that have a power singularity, can be represented
in a neighborhood of the singular point~$x_*$ as
\begin{align}
y\approx A|x_*-x|^{-\beta},\quad \ \beta>0,
\notag
\end{align}
where $A$ is a constant. For these solutions we have $\ds\lim_{x\to x_*}y=\infty$ and $\ds\lim_{x\to x_*}y'_x=\infty$.

The direct application of the standard numerical methods in such problems leads
to certain difficulties because of the singularity in the blow-up solutions
and the unknown (in advance) blow-up point~$x_*$.
Some special methods for solving such problems are described, for example, in \cite{aco2002,mor1979,hir2006,dlam2012}.

Below we propose new methods of numerical integration of such problems.

\section{Problems for first-order equations}

\subsection{Solution method based on a differential transformation}\label{ss:2.1}

The Cauchy problem for the first-order differential equation has the form
\begin{align}
y'_x=f(x,y)\quad (x>x_0);\quad \ y(x_0)=y_0.
\label{eq:02}
\end{align}
In what follows we assume that $f=f(x,y)>0$, $x_0\ge 0$, $y_0>0$, and  $f/y\to\infty$ as $y\to\infty$
(in such problems, blow-up solutions arise when the right-hand side of a nonlinear ODE is
quite rapidly growing as $y\to\infty$).

First, we present the ODE \eqref{eq:02} as a system of equations
\begin{align}
t=f(x,y),\quad \ y'_x=t.
\label{eq:02b}
\end{align}
Then, by applying~\eqref{eq:02b}, we derive a system of equations of the standard form, assuming that
$y=y(t)$ and $x=x(t)$.
By taking the full differential of the first equation in~\eqref{eq:02b} and multiplying
the second one by~$dx$, we get
\begin{equation}
dt=f_x\,dx+f_y\,dy,\quad dy=t\,dx,
\label{eq:02c}
\end{equation}
where $f_x$ and $f_y$ are the respective partial derivatives of $f$.
Eliminating first~$dy$, and then~$dx$ from~\eqref{eq:02c}, we obtain a system of the first-order coupled equations
\begin{equation}
x'_t=\frac 1{f_x+tf_y},\quad \ y'_t=\frac t{f_x+tf_y}\quad (t>t_0),
\label{eq:02d}
\end{equation}
which must be supplemented by the initial conditions
\begin{align}
x(t_0)=x_0,\quad y(t_0)=y_0,\quad t_0=f(x_0,y_0).
\label{eq:02e}
\end{align}

Let $f_x\ge 0$ and $f_y>0$ at $t_0<t<\infty$. Then the Cauchy problem \eqref{eq:02d}--\eqref{eq:02e}
can be integrated numerically, for example, by applying the Runge--Kutta method or other
standard numerical methods (see for example \cite{sch,asch}).
In this case, the difficulties (described in the introduction)  will not occur
since $x'_t$ rapidly tends to zero as $t\to\infty$.
The required blow-up point is determined as $\ds x_*=\lim_{t\to\infty}x(t)$.

\subsection{Examples of test problems and numerical solutions}\label{ss:2.2}

\textit{Example 1}.
Consider the model Cauchy problem for the first-order ODE
\begin{align}
y'_x=y^2\quad (x>0);\quad \ \ y(0)=a,
\label{eq:02f}
\end{align}
where $a>0$. The exact solution of this problem has the form
\begin{align}
y=\frac{a}{1-ax}.
\label{eq:02g}
\end{align}
It has a power-type singularity (a first-order pole) at a point $x_*=1/a$.

By introducing a new variable $t=y'_x$ in~\eqref{eq:02f}, we obtain the following Cauchy problem for
the system of equations:
\begin{equation}
x'_t=\frac{1}{2ty},\ \ y'_t=\frac{1}{2y}\ \ (t>t_0);\quad\ x(t_0)=0,\ y(t_0)=a,\ t_0=a^2,
\label{eq:02h}
\end{equation}
which is a particular case of the problem \eqref{eq:02d}--\eqref{eq:02e} with $f=y^2$, $x_0=0$, and $y_0=a$.
The exact solution of this problem has the form
\begin{align}
x=\frac{1}{a}- \frac{1}{\sqrt{t}},\quad \ y=\sqrt{t}\quad (t\ge a^2).
\label{eq:02i}
\end{align}
It has no singularities; the function $x=x(t)$ increases monotonically for $t>a^2$, tending
to the desired limit value $\ds x_*=\lim_{t\to\infty}x(t)=1/a$, and the function $y=y(t)$ monotonously
increases with increasing~$t$.
The solution~\eqref{eq:02i} for the system~\eqref{eq:02h} is a solution of the original problem~\eqref{eq:02f} in parametric form.

Let $a=1$. Figure~{\sl 1a} shows a comparison of the exact solution~\eqref{eq:02g} of the Cauchy problem for one equation~\eqref{eq:02f}
with the numerical solution of the system of equations~\eqref{eq:02h},
obtained by the classical Runge--Kutta method (with stepsize=0.2).

\begin{figure}
\centering
\vskip-2pc
{\includegraphics[scale=0.34]{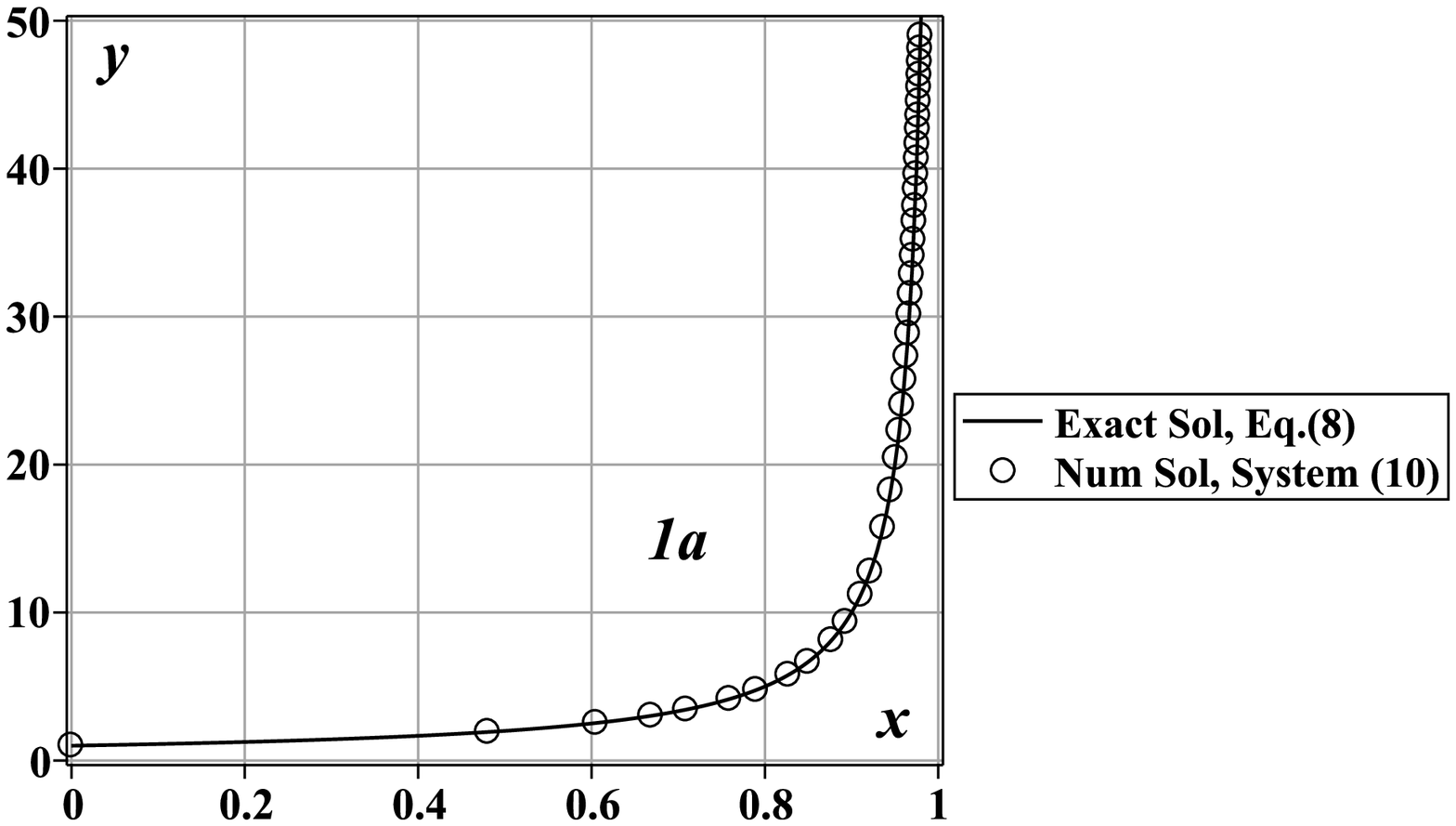}\ \includegraphics[scale=0.34]{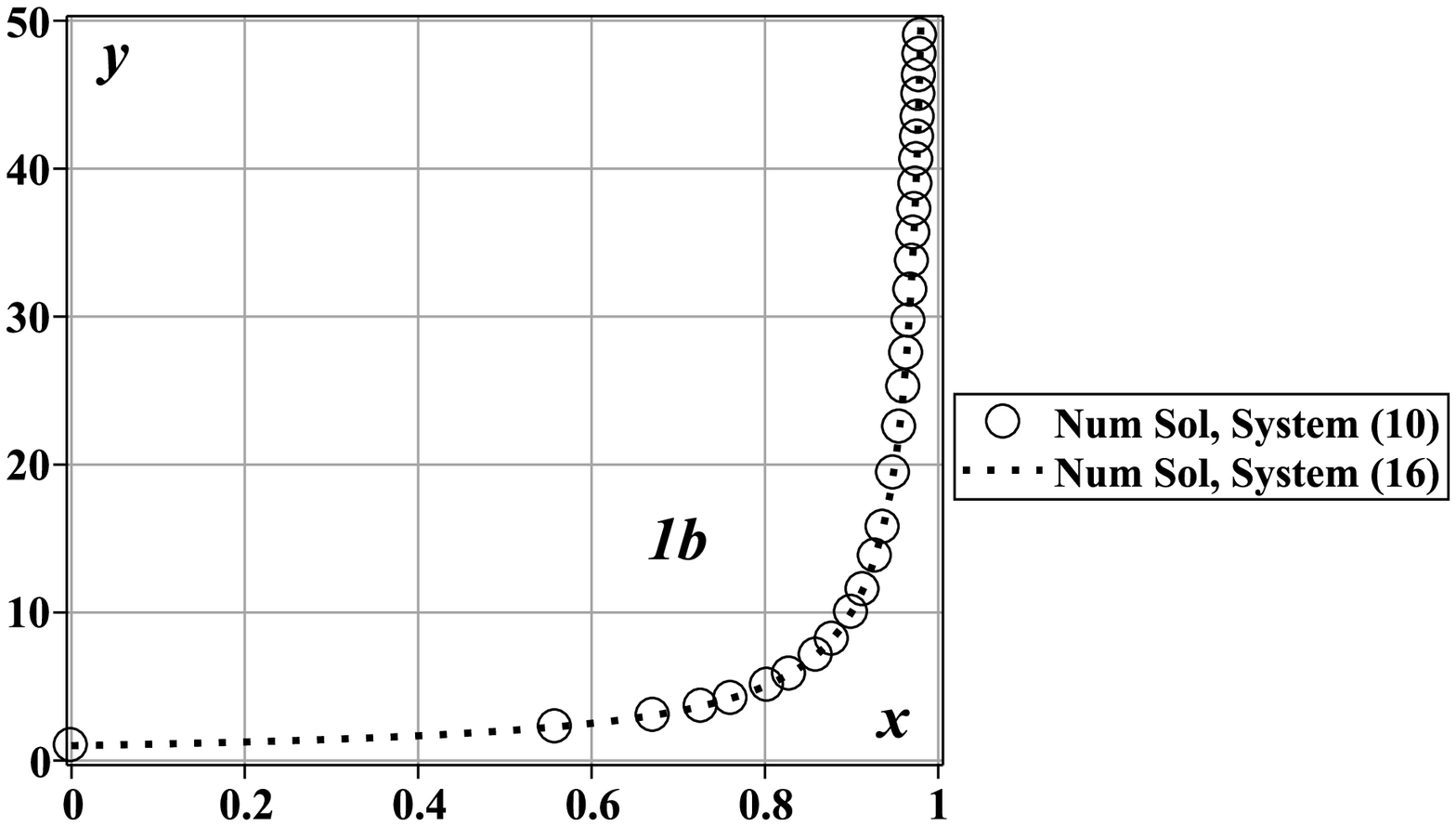}}
\vskip-1pc
\caption{{\sl 1a}\dash Exact solution~\eqref{eq:02g} of the Cauchy problem~\eqref{eq:02f}
and numerical solution of system~\eqref{eq:02h}; {\sl 1b}\dash
numerical solutions of systems~\eqref{eq:02h} and~\eqref{eq:02r} ($a=1$ and $x_*=1$).}
\label{fig:Fig1}
\end{figure}

\subsection{Solution method based on non-local transformations}\label{ss:2.4}

Introducing a new non-local variable according to the formula,
\begin{align}
\xi=\int^x_{x_0}g(x,y)\,dx,\quad \ y=y(x),
\label{*}
\end{align}
leads the Cauchy problem for one equation \eqref{eq:02} to the equivalent problem for the
autonomous system of equations
\begin{equation}
x'_\xi=\frac 1{g(x,y)},\quad \ y'_\xi=\frac {f(x,y)}{g(x,y)}\quad \ \ (\xi>0);\quad \ \ x(0)=x_0,\quad y(0)=y_0.
\label{eq:02p}
\end{equation}
Here, the function $g=g(x,y)$ has to satisfy the following conditions:
\begin{equation}
g>0\ \text{at}\ x\ge x_0, \ y\ge y_0; \quad \ g\to \infty \ \text{as}\ y\to\infty;\quad \ f/g=k\ \text{as}\ y\to\infty,
\label{eq:02q}
\end{equation}
where $k=\text{const}>0$ (and the limiting case  $k=\infty$ is also allowed); otherwise
the function~$g$ can be chosen rather arbitrarily.

From \eqref{*} and the second condition~\eqref{eq:02q} it follows that
$x'_\xi\to 0$ as $\xi\to\infty$. The Cauchy problem~\eqref{eq:02p} can be integrated numerically
applying the Runge--Kutta method or other standard numerical methods.

Let us consider some possible selections of the function~$g$ in the system~\eqref{eq:02p}.

$1^\circ$. We can take $g=\bl(1+|f|^s\br)^{1/s}$ with $s>0$.
In this case, $k=1$ in~\eqref{eq:02q}. For $s=2$,
we get the method of the arc length transformation~\cite{mor1979}.

$2^\circ$.  We can take $g=f/y$ that corresponds to $k=\infty$ in~\eqref{eq:02q}.

\textit{Example 2}.
For the test problem~\eqref{eq:02f}, in which $f=y^2$, we have $g=f/y=y$.
By substituting these functions in~\eqref{eq:02p}, we arrive at the Cauchy problem
\begin{equation}
x'_\xi=\frac{1}{y},\quad \ y'_\xi=y\quad \ \ (\xi>0);\quad \ \
x(0)=0,\quad y(0)=a.
\label{eq:02r}
\end{equation}
The exact solution of this problem is written as follows:
\begin{equation}
x=\frac{1}{a}\bl(1-e^{-\xi}\br),\quad \ y=ae^\xi.
\label{eq:02pe}
\end{equation}
We can see that the unknown quantity $x=x(\xi)$ exponentially tend to the asymptotic values $x=x_*=1/a$ as $\xi\to \infty$.

The numerical solutions of the problems~\eqref{eq:02h} and~\eqref{eq:02r},
obtained by the fourth-order Runge--Kutta method, are presented in Fig.~{\sl 1b} (for $a=1$ and
the same stepsize $=0.2$).
The numerical solutions
are in a good agreement, but the method based on the non-local transformation with $g=t/y$
is more effective than the method based on a differential transformation.

\section{Problems for second-order equations}

\subsection{Solution method based on a differential transformation}\label{ss:3.1}

The Cauchy problem for the second-order differential equation has the form
\begin{align}
y''_{xx}=f(x,y,y'_x)\quad (x>x_0);\quad \ \ y(x_0)=y_0,\quad \ y'_x(x_0)=y_1.\label{eq:03}
\end{align}

Note that the exact solutions of equations of the form~\eqref{eq:03}, which can be used
for test problems with blow-up solutions, can be found,
for example, in \cite{pol2003,kudr}.

Let $f(x,y,u)>0$ if $y>y_0\ge 0$ and $u>y_1\ge 0$, and
the function~$f$ increases quite rapidly as $y\to\infty$
(e.g. if $f$ does not depend on~$y'_x$, then
$\ds\lim_{y\to\infty}f/y=\infty$).

First, we represent ODE~\eqref{eq:03} as an equivalent
system of equations
\begin{align}
y'_x=t,\quad \ y''_{xx}=f(x,y,t).
\label{eq:03b}
\end{align}
where $y=y(x)$ and $t=t(x)$ are unknown functions.
Taking into account~\eqref{eq:03b}, we derive further a standard system of equations, assuming that
$y=y(t)$ and $x=x(t)$.
To do this, we differentiate the first equation in~\eqref{eq:03b} with respect to~$t$.
We have $(y'_x)'_t=1$.
Taking into account the relations $y'_t=tx'_t$ (follows from the first equation of~\eqref{eq:03b})
and $(y'_x)'_t=y''_{xx}/t'_x=x'_ty''_{xx}$, we get further
\begin{equation}
x'_ty''_{xx}=1.
\label{eq:03c}
\end{equation}
If we eliminate the second derivative $y''_{xx}$ by using a second equation of~\eqref{eq:03b},
we obtain the first-order equation
\begin{equation}
x'_t=\frac 1{f(x,y,t)}.
\label{eq:03d}
\end{equation}
Considering further the relation $y'_t=tx'_t$,  we transform~\eqref{eq:03d} to the form
\begin{equation}
y'_t=\frac t{f(x,y,t)}.
\label{eq:03e}
\end{equation}
Equations \eqref{eq:03d} and \eqref{eq:03e} represent a system of coupled first-order equations
with respect to functions $x=x(t)$ and $y=y(t)$. The system \eqref{eq:03d}--\eqref{eq:03e}
should be defined with the initial conditions
\begin{align}
x(t_0)=x_0,\quad y(t_0)=y_0,\quad t_0=y_1.
\label{eq:03f}
\end{align}

The Cauchy problem \eqref{eq:03d}--\eqref{eq:03f} can be integrated numerically
applying the standard numerical methods \cite{sch,asch}, without fear of blow-up solutions.


\begin{remark}
Systems of equations
\eqref{eq:02b} and \eqref{eq:03b} are particular cases of parametrically defined
nonlinear differential equations, which are considered in \cite{pol2016,pol2017}.
In \cite{pol2017}, the general solutions of several parametrically defined ODEs were obtained
via differential transformations, based on introducing a new independent
variable~$t=y'_x$.
\end{remark}

\subsection{Examples of test problems and numerical solutions}\label{ss:3.2}

\begin{figure}
\centering
\vskip-2pc
{\includegraphics[scale=0.34]{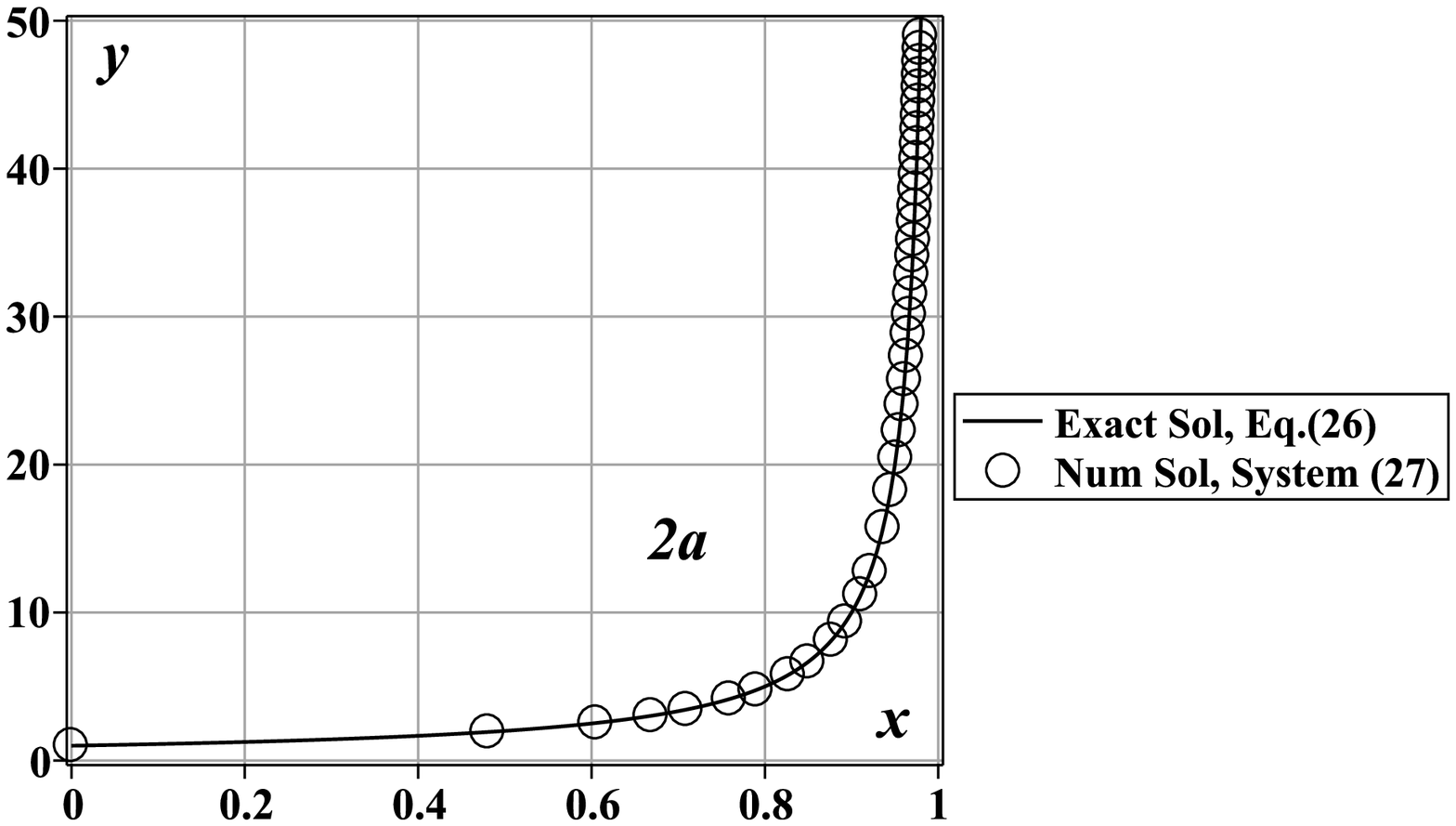}\ \includegraphics[scale=0.34]{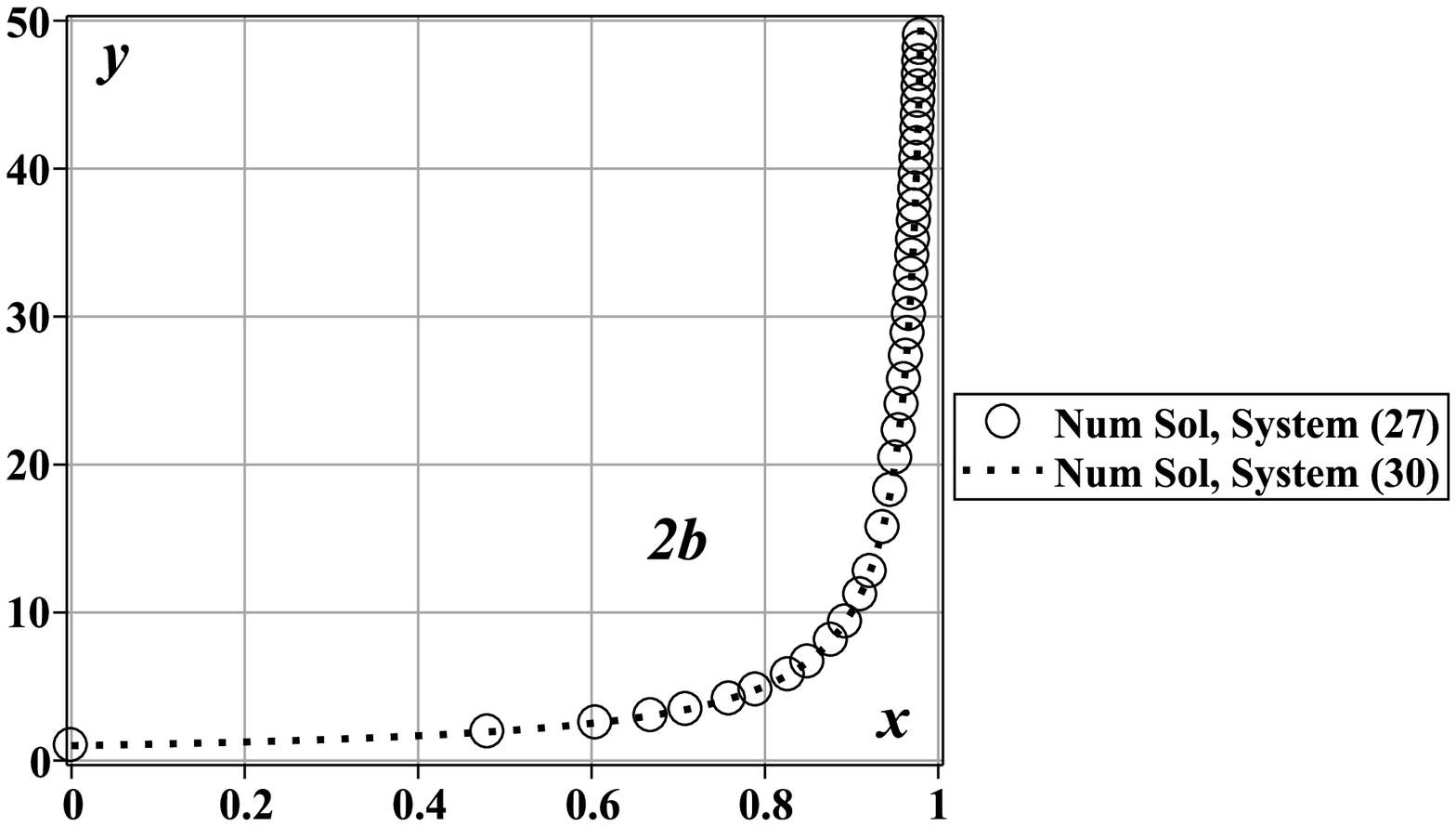}}
\vskip-1pc
\caption{{\sl 2a}\dash Exact solution~\eqref{eq:02g} of the Cauchy problem~\eqref{eq:04a}
and numerical solution of system~\eqref{eq:04b};
{\sl 2b}\dash numerical solutions of systems~\eqref{eq:04b} and~\eqref{eq:xx} ($a=1$ and $x_*=1$).}
\label{fig:Fig2}
\end{figure}

\textit{Example 3}.
Let us consider Cauchy problem
\begin{align}
y''_{xx}=2y^3\quad \ (x>0),\quad \ y(0)=a,\quad \ y'_x(0)=a^2.
\label{eq:04a}
\end{align}
The exact solution of this problem is defined by the formula~\eqref{eq:02g}.
\goodbreak

Introducing a new variable $t=y'_x$ in~\eqref{eq:04a}, we transform~\eqref{eq:04a} to the Cauchy problem
for the system of the first-order ODEs
\begin{equation}
x'_t=\tfrac 12y^{-3}, \ \ \ y'_t=\tfrac12 ty^{-3} \ \ \ (t>t_0);\quad \
x(t_0)=0, \ \ \ y(t_0)=a, \ \ \ t_0=a^2,
\label{eq:04b}
\end{equation}
which is a particular case of the problem \eqref{eq:03d}--\eqref{eq:03f} with $f=y^2$, $x_0=0$, and $y_0=a$.
The exact solution of the problem~\eqref{eq:04b} is given by the formulas~\eqref{eq:02i}.

Figure~{\sl 2a} shows a comparison of the exact solution~\eqref{eq:02g}  of the Cauchy problem for one equation~\eqref{eq:04a}
with the numerical solution of the system of equations~\eqref{eq:04b},
obtained by the fourth-order Runge--Kutta method (we have a good coincidence).

\subsection{Solution method based on non-local transformations}\label{ss:3.3}

First, equation~\eqref{eq:03} can be represented as a system of two equations
$$
y'_x=t,\quad \ t'_x=f(x,y,t),
$$

\noindent
and then we introduce the non-local variable by the formula
\begin{align}
\xi=\int^x_{x_0}g(x,y,t)\,dx,\quad \ y=y(x),\quad t=t(x).
\label{eq:05}
\end{align}

\noindent
As a result, the Cauchy problem \eqref{eq:03} can be transformed to
the following problem for an autonomous system of three equations:
\begin{equation}
\begin{gathered}
x'_\xi=\frac 1{g(x,y,t)},\quad \ y'_\xi=\frac t{g(x,y,t)},\quad \ t'_\xi=\frac {f(x,y,t)}{g(x,y,t)}\quad \ \ (\xi>0);\\
x(0)=x_0,\quad \ y(0)=y_0,\quad \ t(0)=y_1.
\end{gathered}
\label{eq:06}
\end{equation}
For a suitable choice of the function $g=g(x,y,t)$ (not very restrictive conditions of the form~\eqref{eq:02q} must be
imposed on it), the Cauchy problem \eqref{eq:06} can be numerically integrated applying the
standard numerical methods \cite{sch,asch}.

Let us consider some possible selections of the function~$g$ in system~\eqref{eq:06}.
\goodbreak

$1^\circ$. We can take $g=\bl(1+|t|^s+|f|^s\br)^{1/s}$ with $s>0$. The case $s=2$
corresponds to the method of the arc length transformation~\cite{mor1979}.

$2^\circ$. Also, we can take $g=f/y$, $g=f/t$, or $g=t/y$.

\textit{Example 4}. For the test problem~\eqref{eq:04a}, in which $f=2y^3$,
we put $g=t/y$. By substituting these functions in~\eqref{eq:06}, we arrive
at the Cauchy problem
\begin{equation}
x'_\xi=y/t,\ \ y'_\xi=y,\ \ t'_\xi=2y^4/t\ \ (\xi>0);\ \ \ x(0)=0,\ \ y(0)=a,\ \ t(0)=a^2.
\label{eq:xx}
\end{equation}
The exact solution of this problem is written as follows:
\begin{equation}
x=a^{-1}\bl(1-e^{-\xi}\br),\quad \ y=ae^\xi,\quad \ t=a^2e^{2\xi}.
\label{eq:xy}
\end{equation}
\noindent
We can see that the unknown quantity $x=x(\xi)$ exponentially tend to the asymptotic values $x=x_*=1/a$ as $\xi\to \infty$.

The numerical solutions of the problems~\eqref{eq:04b} and~\eqref{eq:xx},
obtained by the fourth-order Runge--Kutta method, are presented in Fig.~{\sl 2b}
(for $a=1$ and the same stepsize $=0.2$).
The numerical solutions
are in a good agreement, but the method based on the non-local transformation with $g=t/y$
is more effective than the method based on a differential transformation.

\begin{remark}
The method described in Section~\ref{ss:3.3} can be generalized to nonlinear ODEs of arbitrary
order and systems of coupled ODEs.
\end{remark}

\end{document}